\documentclass[letterpaper, 10 pt, conference]{ieeeconf}
\usepackage{enumerate}

\usepackage{cite}
\usepackage{amsmath,amssymb,amsfonts}
\usepackage{algorithmic}
\usepackage{graphicx}
\usepackage{textcomp}
\usepackage{fancyhdr}

\usepackage{siunitx}
\usepackage[short]{optidef}
\usepackage[pdfauthor={Jakob Harzer},
pdftitle={Efficient Optimal Control for Highly Oscillatory Systems},
pdfstartview=XYZ,
bookmarks=true,
colorlinks=true,
linkcolor=black,
urlcolor=blue,
citecolor=blue,
pdftex,
bookmarks=true,
linktocpage=true, 
hyperindex=true,
]{hyperref}

\newcommand{\R}{{\mathbb{R}}}

\newcommand{\dd}[1]{\ \mathrm{d}#1}

\DeclareMathSymbol{\shminus}{\mathbin}{AMSa}{"39}

\newcommand{\smallPlus}{{\scalebox{.7}{$+$}}}
\newcommand{\smallMinus}{{\scalebox{0.6}[0.7]{$-$}}}

\newcommand{\xminus}{x^\smallMinus}
\newcommand{\xplus}{x^\smallPlus}

\def\BibTeX{{\rm B\kern-.05em{\sc i\kern-.025em b}\kern-.08em T\kern-.1667em\lower.7ex\hbox{E}\kern-.125emX}}

\fancypagestyle{FirstPage}{
\lfoot{
    \noindent\fbox{%
    \parbox{\textwidth}{%
    \footnotesize\copyright2022 IEEE.  Personal use of this material is permitted.  Permission from IEEE must be obtained for all other uses, in any current or future media, including reprinting/republishing this material for advertising or promotional purposes, creating new collective works, for resale or redistribution to servers or lists, or reuse of any copyrighted component of this work in other works.
    }%
}    
} 
}

\begin{document}
\title{Efficient Numerical Optimal Control for Highly Oscillatory Systems}
\author{Jakob Harzer, Jochem De Schutter, Moritz Diehl
\thanks{This research was supported by DFG via Research Unit FOR 2401 and project 424107692 and by the EU via ELO-X 953348.
Jakob Harzer and Jochem De Schutter are with the Department of Microsystems Engineering (IMTEK) and Moritz Diehl is with the Department of Microsystems Engineering (IMTEK) and Department of Mathematics, University of Freiburg, Georges-Koehler-Allee 102, 79110 Freiburg, Germany.
Emails: \textbf{harzer.jakob@gmail.com}, jochem.de.schutter@imtek.uni-freiburg.de, moritz.diehl@imtek.uni-freiburg.de
}}
\maketitle
\begin{abstract}
We present an efficient transcription method for highly oscillatory optimal control problems. For these problems, the optimal state trajectory consists of fast oscillations that change slowly over the time horizon. 
Out of a large number of oscillations, we only simulate a subset to approximate the slow change by constructing a semi-explicit differential-algebraic equation that can be integrated with integration steps much larger than one period. 
For the solution of optimal control problems with direct methods, we provide a way to parametrize and regularize the controls.
Finally, we utilize the method to find a fuel-optimal orbit transfer of a low-thrust satellite. 
Using the novel method, we reduce the size of the resulting nonlinear program by more than one order of magnitude. 
\end{abstract}
\thispagestyle{FirstPage}
\section{Introduction}
Oscillating motion appears in a wide range of physical systems, from fast oscillations in crystal oscillators, over the forced vibrations of an audio speaker, up to the very slow tides. 
Although we might not expect it to, the oscillating motion might not be periodic but also include a slow change, that, compared to the fast oscillations, is barely noticeable. 
To notice this change, one would have to observe the motion over a long horizon and possibly a large number of oscillations. 

We will call this behavior `highly oscillatory' and we define `highly oscillatory systems' following \cite{petzold_jay_yen_1997} as  ``systems whose solutions may be oscillatory in the sense that there is a fast solution which varies regularly about a slow solution'' under the condition that the timescale of the fast solution is much shorter than the interval of interest. 
In the same fashion, we define the property  `highly oscillatory' for a trajectory of a controlled system.

When simulating highly oscillatory systems, a sufficiently small integration step size is needed to simulate the system accurately over a single cycle. 
But since we are interested in the long-term behavior over a large number of cycles, the number of integration steps required becomes excessively large.
We can transfer this definition also to optimal control problems (OCP).
A continuous-time OCP is now called `highly oscillatory' if the optimal state and control trajectories
are expected to be highly oscillatory.

When a highly oscillatory OCP is solved using direct multiple shooting, a large number of steps and controls corresponds to a large number of variables and constraints in the resulting nonlinear program (NLP), which in consequence requires huge computational effort to solve.

The question we want to answer in this paper is:
\begin{quote}
    ``How can we exploit the fact that the expected optimal trajectory of a system in an optimal control problem is highly oscillatory in order to speed up the optimization process?''
\end{quote}

\begin{figure}%
    \centering %
    \includegraphics*[width=\linewidth]{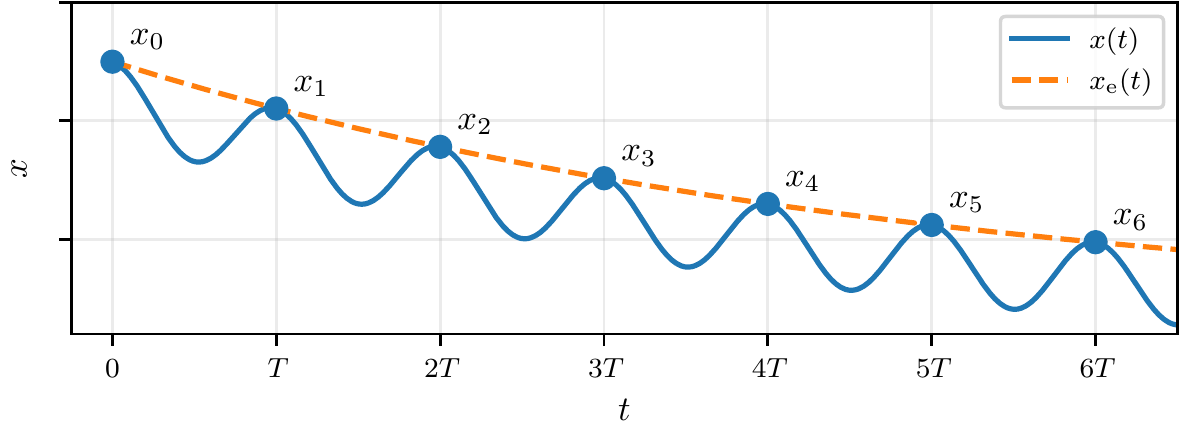}%
    \vspace{-7pt}%
    \caption{Highly oscillatory trajectory and its envelope}
    \label{fig:ExampleHOsystem} %
\end{figure}

Consider the simulation of a system \(f: \R^{n_x} \times \R^{n_u} \rightarrow \R^{n_x}\). 
Suppose \(x(t)  \in \R^{n_x},\ t \in [0,t_\mathrm{f}]\) is the highly oscillatory solution of the initial value problem \begin{subequations}
    \begin{align}
        \dot{x}(t) &= f(x(t),u(t)) \label{eq:IVPrhs} \\
        x(0) &= x_0
    \end{align}
    \label{eq:IVP}%
\end{subequations} where \(x_0 \in \R^{n_x}\), as visualized in Fig. \ref{fig:ExampleHOsystem}. 
Let us assume for a moment that the controls are fixed, that we know the period of the oscillations \(T\), and that we have a one-period map \(\Psi: \R^{n_x} \rightarrow \R^{n_x}\) that computes the state \(x_{k+1} = \Psi(x_k)\) after one period from a given \(x_k\).
Since the cycles only change slowly, the map \(\Psi\) is a near-to-identity map \cite{Calvo2012}. 
Starting from \(x_0\) we are able to generate a series of points \(x_1,x_2,\dots\) that samples the trajectory periodically. 
In the manner of \cite{Petzold}, a trajectory that smoothly interpolates this slowly changing series is called the `envelope' \(x_\mathrm{e}(t)\).
We assume that \(x_\mathrm{e}(t)\) is `sufficiently smooth' in the sense that its derivatives in \(t\) exist and are bounded, and that it stems from some envelope dynamics%
\begin{align}%
    \dot{x}_\mathrm{e}(t) = f_e(x_\mathrm{e}(t)) \label{eq:envelopeDynamics}
\end{align}%
whose solution \(x_\mathrm{e}(t)\) is equal to the highly oscillatory solution \(x(t)\) at multiples of the period \(kT, k \in \mathbb{N}\). 
More details can be found in \cite{Calvo2012}. 
Since these envelope dynamics are assumed to generate a slowly changing trajectory, we can integrate them with stepsizes much bigger than a single cycle duration, while still retaining a decent accuracy. 

\subsection{Contribution}

In this work, we construct a differential-algebraic equation (DAE) that approximates the envelope dynamics of the highly oscillatory system.
Instead of simulating all cycles of an oscillating trajectory, we only simulate a few to gather information about the slow change that occurs over the time horizon. 
To account for unknown and possibly changing cycle durations \(T\) of the oscillations, we will scale the independent variable \(t\) and use phase conditions to define a cycle.
We provide a way to parametrize the controls and use the method to efficiently solve highly oscillatory OCPs using direct methods. 
This method can significantly reduce the number of variables in the resulting NLP, but requires careful regularization of the state trajectories.

This simulation approach is closely related to the \textit{stroboscopic method for highly oscillatory problems} \cite{Calvo2012} where the envelope dynamics are also modeled with a differential equation and then numerically approximated. 
Quite similar as well are the \textit{multirevolution methods}\cite{Graf1972} and the \textit{envelope following methods} \cite{Petzold,petzold_jay_yen_1997,Graf1972,Gallivan1980} that instead of approximating the continuous system \eqref{eq:envelopeDynamics}, model the envelope dynamics using a difference equation. 

We propose an envelope-following formalism specifically tailored for direct optimal control.
Whereas averaging methods have been used for highly oscillatory optimal control before, such as in \cite{Yang2009}, to the authors' knowledge, ideas from envelope following methods, on the other hand, have not been used for optimal control yet.

\subsection{Outline}
After introducing a guiding example, we start with the definition of a single cycle of the oscillations in Section \ref{sec:CycleDefinition}. In Section \ref{sec:NCyclesOCP} we construct an \(N\)-cycle OCP, that we then simplify by approximating the envelope dynamics using a DAE in Section \ref{sec:EnvelopeDAE}. 
We discuss the integration of the DAE and problems of this approach in Section \ref{sec:MacroIntegration}.
In Section \ref{sec:Satellite} we discuss a numerical example on the model of a low-thrust satellite, and conclude the paper with Section \ref{sec:Conclusion}.

All the coding of this work is done in Python 3.9, we use CasADi \cite{Andersson2019} to set up the optimization problems and solve the nonlinear programs using IPOPT \cite{IPOPT} equipped with the linear solver MA27\cite{HSLMA27}.

\subsection{Guiding Example}

As a guiding example, we consider the nonlinear predator-prey system, also known as the Lotka-Volterra equations. The system models the interaction between a prey population of size \(r\) and a predator population of size \(s\).
We modify the original system by adding a control \(u \in \R \) to influence the growth of the predator species. The model is then given by
\begin{align}
    x := \begin{bmatrix}r \\ s \end{bmatrix} \in \R^2,  && f(x,u) := \begin{bmatrix} \frac{2}{3} r - \frac{4}{3} r s \\ r s - s \end{bmatrix} + \begin{bmatrix} 0 \\ u \end{bmatrix}. \label{eq:PredatorPreyModelEquations}
\end{align}
For visualization, Fig. \ref{fig:PhaseConditions1} shows a trajectory  from the initial state \(\bar{x}_0 = [1,1.5]^\top\) with constant control \(u(t) = 0.03\).


\section{Definition of a single cycle}
\label{sec:CycleDefinition}%
\color{black} %
We discretize the oscillatory trajectory into individual cycles by introducing two affine phase conditions 
\begin{subequations}
    \begin{align}
        q^\top \xminus = b^\smallMinus \label{eq:LinearPhase1} \\
        q^\top \xplus = b^\smallPlus \label{eq:LinearPhase2}
    \end{align} %
\end{subequations}
for the starting point \(x^\smallMinus \in \R^{n_x}\) and end point \(x^\smallPlus \in \R^{n_x}\) of a cycle, respectively, that are defined by  a normalized vector \(q \in \R^{n_x}\) that satisfies \(\|q\| = 1\), and scalars \(b^\smallMinus,b^\smallPlus \in \R\).
These phase conditions have to be chosen depending on the system, and usually involve prior knowledge of the expected oscillatory trajectories. 
We distinguish two cases:
\begin{enumerate}[(A)]
    \item The trajectory \(x(t)\) is discretized by giving a condition on an oscillatory state. For these types of systems we can choose \(b^\smallMinus = b^\smallPlus\), which is equivalent to defining a linear Poincaré section. 
    For example, for the predator-prey system \eqref{eq:PredatorPreyModelEquations} we will specify that ``for the starting and end point of a cycle the first state is 1'' as the phase conditions\begin{subequations}%
        \begin{align}%
            q^\top & \xminus - 1 = 0  \\
            q^\top & \xplus - 1 = 0,
        \end{align} \label{eq:PhaseConditionsPredPrey} %
    \end{subequations} %
    for which we choose \(q = [1,\ 0]^\top\) and \(b^\smallPlus = b^\smallMinus = 1\). This condition is visualized in Fig. \ref{fig:PhaseConditions1}.
    \begin{figure}
        \vspace{4pt}
        \centering%
        \includegraphics[width=\linewidth]{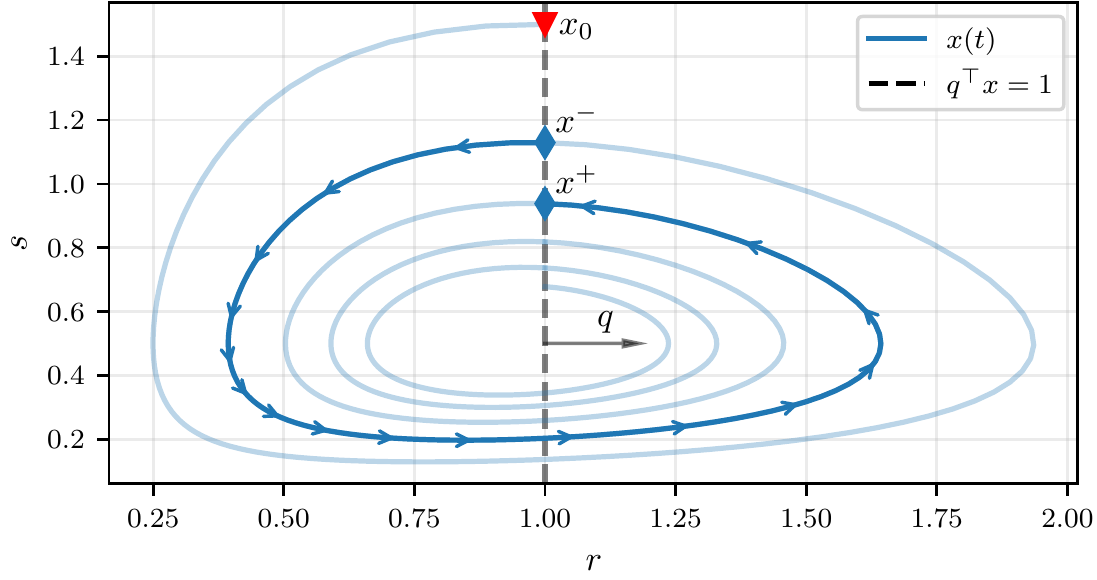}%
        \vspace{-7pt}%
        \caption{Phase conditions \eqref{eq:PhaseConditionsPredPrey} of type (A) discretize the state trajectory of the Predator-Prey system into individual cycles.} %
        \label{fig:PhaseConditions1} %
    \end{figure} %
    \item The trajectory \(x(t)\) is discretized by giving a condition on a \textit{non-oscillatory state}, typically a `phase state' such as a clockstate or a state for a rotation angle. 
    This is required if this state influences the dynamics strongly and periodically over a single cycle, cf. \cite{Gallivan1980,Calvo2012}.
    This approach is used later in the example in Section \ref{sec:Satellite}. 
\end{enumerate}
Similar phase conditions have also been used by others to define a cycle of the oscillations \cite{Ting2008,Linaro2020,Maffezzoni2009,Petzold,Gallivan1980,Gear1982}.


To later integrate over \(N\) cycles and also to handle changing cycle durations, we introduce the numerical time \(\tau \in [0, N]\), similar to \cite{Petzold,petzold_jay_yen_1997}.
The time \(t \in [0,t_\mathrm{f}]\) will now be called the `physical time' and can be included as an additional state in the model. 
We define, by a slight abuse of notation, \(x(\tau) := x(t(\tau))\).
The relationship between \(t\) and \(\tau\) is given by 
\begin{align}
    \frac{\dd t(\tau)}{\dd \tau} = T(\tau), && T(\tau) \coloneqq T_k,\ \tau \in [k,k+1),
\end{align} %
where \(T_k\) is the duration of the \(k\)-th cycle.
This way, expressed in numerical time, the same control grid can be applied to each cycle.

Let us parametrize the controls over one cycle by a matrix
\begin{align}
    U_k = [u_{k,0},u_{k,1},\dots,u_{k,N_\mathrm{ctr}-1}] \in \R^{n_u \times N_\mathrm{ctr}}
\end{align}
and apply 
\(u(\tau) = \tilde{u}(\tau; U_k)\) for \(\tau \in [k,k+1)\).  
For example, as depicted in Fig. \ref{fig:Full_OCP}, for the cycle \(k\), this could be a piecewise constant control parameterization 
\begin{align}
\tilde{u}(\tau; U_k) = u_{k,m},\ \text{ for }\ \tau\ \text{mod}\ k \in \left[\tfrac{m}{N_\mathrm{ctr}},\tfrac{m + 1}{N_\mathrm{ctr}}\right),
\end{align}
and for \(m = 0,..,N_\mathrm{ctr}-1\). A single cycle is divided into \(N_\mathrm{ctr}\) equidistant intervals and the \(m\)-th control is applied over the \(m\)-th interval.
The scaled dynamics for a cycle with control parameterization \(U_k\) now read 
\begin{align}
    \frac{\dd x}{\dd \tau}(\tau) = T(\tau) f(x(\tau),\tilde{u}(\tau; U_k)), \qquad  \tau \in [k,k+1).
\end{align} 
Let \(F\) be a numerical integration routine that simulates these dynamics over the numerical cycle duration \(\Delta\tau = 1\). 
In the fashion of \cite{Calvo2012}, we call this integrator the `micro-integrator'.

To numerically simulate a single cycle with a given control \(U_k\) that starts at a point \(x_k\), we solve\begin{subequations}
    \begin{align} 
        0 &= Q^\top(\xminus_k - x_k) \\
        0 &= C(z_k, U_k),%
    \end{align}    %
\end{subequations}%
where the matrix \(Q\) is defined such that \(\bar{Q} = [q | Q]\) is an orthonormal basis of the state space.

Here, the algebraic equations 
\begin{align}
    C(z,U) \coloneqq \begin{cases}
        0 &= q^\top\xminus - b^\smallMinus  \\
        0 &= \xplus - F(\xminus,U,T) \\
        0 &= q^\top\xplus - b^\smallPlus 
\end{cases} \qquad \in \R^{n_x + 2}
    \label{eq:CycleConditions} %
\end{align}
define one cycle of duration \(T\) with starting point \(\xminus\) and end point \(\xplus\) where the variables are collected in \(z \coloneqq (\xminus,\xplus,T)\).
These equations can be solved numerically using for example a Newton-type root-finding method.
By construction system \eqref{eq:CycleConditions} has multiple periodic solution, hence a sufficiently good initial guess is required to simulate exactly one cycle.



%


\section{\(N\)-Cycle Optimal Control Problem}
\label{sec:NCyclesOCP}
\begin{figure}
    \centering 
    \vspace{4pt}
    \includegraphics[width=\linewidth]{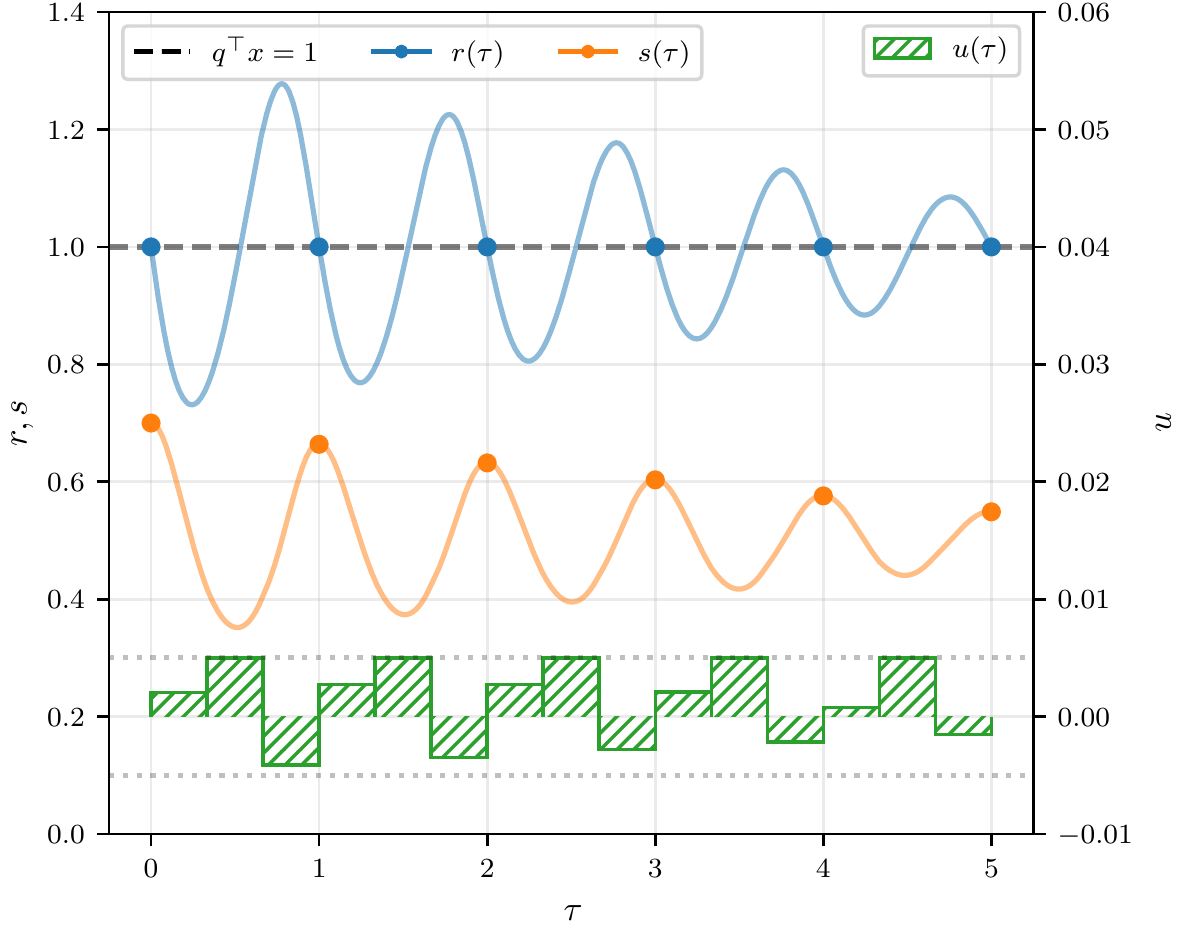} %
    \vspace{-17pt}%
\caption{Exemplary optimal solution of a 5-cycles OCP of the predator-prey system \eqref{eq:PredatorPreyModelEquations}.}%
    \label{fig:Full_OCP} %
\end{figure}

By discretizing the state and control trajectories as explained above, we can formulate an optimal control problem over a number of \(N\) cycles as
\begin{mini!}
    {w}{\sum_{i=0}^{N-1} L(z_k,U_k) + E(\xplus{N-1})}
    {\label{eq:FullOCP}}{}
    \addConstraint{0}{= Q^\top (\xminus_0  - x_0} \label{eq:FullOCP_Initial}) 
    \addConstraint{0}{= C(z_k,U_k) \quad & k = 0,..,N-1 \label{eq:FullOCP_CycleCond}}
    \addConstraint{0}{\leq H(z_k,U_k) \quad & k = 0,..,N-1 \label{eq:FullOCP_Inequalities}}
    \addConstraint{0}{= Q^\top(\xminus_{k} - \xplus_{k-1}) \quad& k = 1,..,N-1 \label{eq:FullOCP_ConnCond}}
\end{mini!}
where \(w = (z_0,\dots,z_{N-1},U_0\dots,U_{N-1})\), \(L\) is a cycle cost and \(E\) a terminal cost.
Essentially, we simulate \(N\) cycles with \eqref{eq:FullOCP_CycleCond} that are constrained to start at the initial point by \eqref{eq:FullOCP_Initial} and are connected by \eqref{eq:FullOCP_ConnCond}.
Equation \eqref{eq:FullOCP_Inequalities} enforces state and control constraints.
Fig. \ref{fig:Full_OCP} visualizes an example of such an \(N\)-cycle OCP. 
As mentioned in the introduction, solving such a problem becomes increasingly costly for larger \(N\).

An alternative control parameterization that accommodates the highly oscillatory nature of the system, derives the cycle controls from a low-order polynomial
\begin{align}
    U(\tau; U_c) = \sum_{i=1}^{d_u} U_{\mathrm{c},i} \ell_i(\tau)
\end{align}
that we evaluate at the cycle \(k\) as \(U_k := U(k+\frac{1}{2}; U_c)\). We construct this polynomial from \(d_u\) construction points  where \(\ell_i\) are the orthogonal basis polynomials built from the construction times \(\tau_{\mathrm{c},1},\dots,\tau_{\mathrm{c},d_u}\). 
The construction controls \(U_c : = (U_{\mathrm{c},1},\dots,U_{\mathrm{c},d_u}) \in \R^{d_u \times n_u \times N_\mathrm{ctr}}\) are included as optimization variables. 
For a  constant periodic control parameterization \(U(\tau) = U_\mathrm{const}\), \(d_u = 1\), i.e. the same periodic control is then applied in each cycle.

\section{Envelope Slope Approximation with a DAE}
\label{sec:EnvelopeDAE}
\begin{figure}
    \centering
    \vspace{4pt}
    \includegraphics[width=\linewidth]{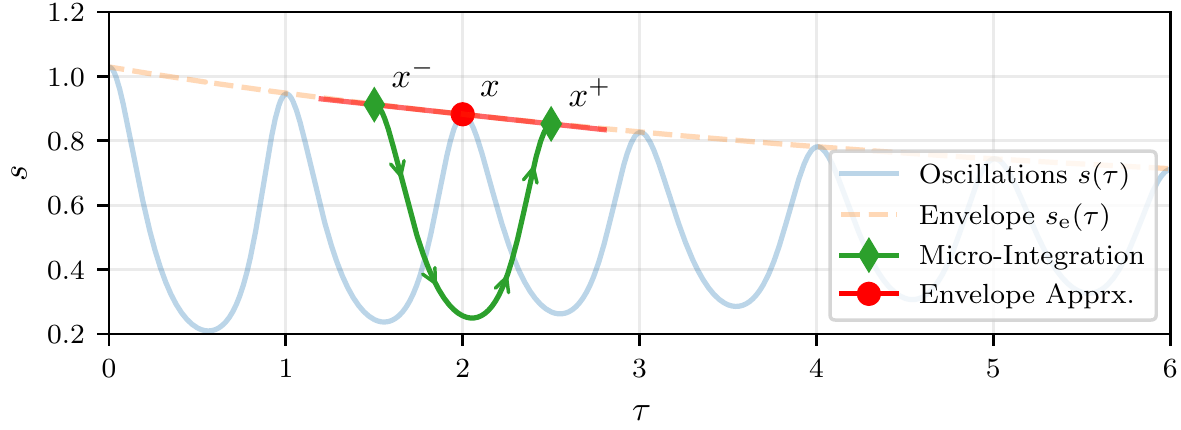}
    \vspace{-17pt}%
    \caption{The 'Central Difference' DAE approximates the dynamics of the envelope of the state \(s\) of the predator-prey model.}
    \label{fig:DAE_CD_sketch}
\end{figure}

The core idea of this work is that instead of integrating the oscillating ODE \eqref{eq:IVPrhs} with a small stepsize, we integrate the DAE
\begin{subequations}
    \begin{align}
        \frac{\dd x}{\dd \tau} &= f_\mathrm{e}(x(\tau),z(\tau),U(\tau)) \\
            0 & = g_\mathrm{e}(x(\tau),z(\tau),U(\tau))
        \end{align}
\end{subequations}
that approximates the envelope dynamics with a large stepsize.
In this paper, we present the `central differences' DAE that approximates the dynamics of the envelope at the point \(x\) as
\begin{subequations}
    \begin{align}
        f_\mathrm{e}(x,z) :=  \xplus - \xminus && \text{(envelope approx.)}  
    \end{align}
    \text{with the algebraic equations \(g_\mathrm{e}(x,z)\) given by}
    \begin{align}
        0 &= Q^\top \left(\frac{\xplus + \xminus}{2} - x\right) & \text{(connecting cond.)} \label{eq:DAE_connCond} \\
        0 &= C(z,U) & \text{(cycle cond.).}  \label{eq:DAE_cyclCond}
    \end{align} %
    \label{eq:DAE_CD} %
\end{subequations} %
Essentially, at a point \(x\), we simulate one cycle around it to estimate the slow change of the highly oscillatory trajectory. 
Algebraic equation \eqref{eq:DAE_cyclCond} corresponds to the micro-integration of a single cycle. 
The approximation of the envelope dynamics is inspired by a central difference scheme
and is sketched in Fig. \ref{fig:DAE_CD_sketch}.
The differential state \(x\) is constrained to `lie in the middle' of the algebraic states \(\xminus\) and \(\xplus\) by the connecting condition \eqref{eq:DAE_connCond}.
Multiplying by \(Q^\top\) in \eqref{eq:DAE_connCond} projects the condition into the null space of \(q\). 
This prevents the DAE from being overdetermined since the direction \(q\) is already constrained by the phase conditions inside \eqref{eq:DAE_cyclCond}.

In a slightly different formulation, a similar method was already used in \cite{Ting2008} to simulate electrical oscillators where the authors employ forward-difference and backward-difference schemes and also suggest using a central-difference scheme. 
Also noteworthy is the approach by \cite{Calvo2012} where the authors approximate the envelope dynamics using the previous and the next cyclic point, that in consequence requires the simulation of two cycles.

\section{Macro-Integration of the Envelope DAE}
\label{sec:MacroIntegration}

For the simulation of the oscillatory system over a timescale much larger than the duration of a single cycle, instead of simulating a large number of \(N\) cycles in detail, we now just integrate the DAE whose solution approximates the envelope of the oscillations. Again, in the manner of \cite{Calvo2012}, this integration procedure will be called the `macro-integration'. 
Since we now use the numerical time \(\tau\) instead of the physical time, integrating the DAE with a stepsize of \(N\) simulates the highly oscillatory dynamics over \(N\) cycles. 
In practice, for the macro-integration of the DAE we use implicit RK methods.

Similar to \cite{Calvo2010}, we can list the following three sources of error for the overall simulation procedure:
\begin{enumerate}[1.]
    \item Errors of the micro-integration over a single cycle.
    \item Errors of the central difference approximation of the envelope dynamics.
    \item Errors of the macro-integration of the constructed envelope DAE.
\end{enumerate}
To analyze these errors analytically is outside the scope of this paper.

\subsection{Simulation Error and Regularization}
\label{sec:IntegrationErrorandRegularization}

The simulation method presented here is tailored to highly oscillatory problems where we expect that the trajectories are composed of fast oscillations that change very little over the time horizon. 
Thus, we expect the cycle duration to only drift slowly, the shape of the cycles to remain similar and the envelope to be sufficiently smooth.  
If this is not the case, we introduce a large simulation error contribution from error source 2.

Dependent on the cost functions \(L\) and \(E\), the optimizer might prefer trajectories with large values in the higher-order derivatives.
This is especially a problem for the presented method that becomes inaccurate if the envelope is not sufficiently smooth. One way to ensure this is to add additional constraints or regularization terms to the objective function that penalize the derivatives of the envelope or control trajectories in the macro-integration.

\subsection{Regularization Experiments}
\label{sec:Experiments}

Suppose we want to drive the predator prey system \eqref{eq:PredatorPreyModelEquations} from an initial state \(\bar{x}_0 = [1,2] \) into a target state \(\bar{x}_\mathrm{target} = [1,0.7]\) within \(N=20\) cycles that are characterized by the phase conditions \eqref{eq:LinearPhase1}.
The controls are very small and bounded by
\begin{align}
    -0.005 \leq u(\tau) \leq 0.005, \forall \tau,
\end{align} such that even with an optimal control, the system cannot be driven to the target point.
Instead of solving the \(N\)-cycle OCP \eqref{eq:FullOCP}, we instead solve an OCP
\begin{mini*}
    {\substack{x_\mathrm{e}(\cdot),z_\mathrm{e}(\cdot),\\U_c}}{\|x_\mathrm{e}(N) - \bar{x}_N\|^2_2 + \alpha \int_0^N \left\|
        \frac{\dd^3 x}{\dd \tau^3}(\tau) \right\|^2_2 \dd\tau}
    {\label{eq:PredatorPreyOCP}}{}
    \addConstraint{0}{= x_\mathrm{e}(0) - \bar{x}_0} 
    \addConstraint{\dot{x_\mathrm{e}}}{= f_\mathrm{e}(x_\mathrm{e}(\tau),z_\mathrm{e}(\tau),U(\tau; U_c)), &\qquad \tau \in [0,N]}
    \addConstraint{0}{= g_\mathrm{e}(x_\mathrm{e}(\tau),z_\mathrm{e}(\tau),U(\tau; U_c)), &\qquad \tau \in [0,N]}
    \addConstraint{0}{\leq h(x_\mathrm{e}(\tau),z_\mathrm{e}(\tau),U(\tau; U_c)), &\qquad \tau \in [0,N]}
\end{mini*} that uses the DAE dynamics.
We search for an optimal polynomial control \(U(\tau; U_c)\) of degree 2, i.e. \(d_u=3\), with \(N_\mathrm{ctrl} = 3\) controls per cycle. 
In each micro-integration, we chain 10 RK4 integration steps to simulate the model dynamics over each of the 3 control intervals, i.e. we perform 30 RK4 steps per cycle.
To keep the envelope sufficiently smooth, we regularize the third derivative of the envelope with a small positive weight \(\alpha > 0\).
We discretize the OCP using a single interval of a 3-stage Gauss-Legendre collocation scheme of order 6 and solve the NLP using the NLP solver IPOPT\cite{IPOPT}.
\begin{figure}
    \vspace{4pt}
    \centering%
    \includegraphics*[width = \linewidth]{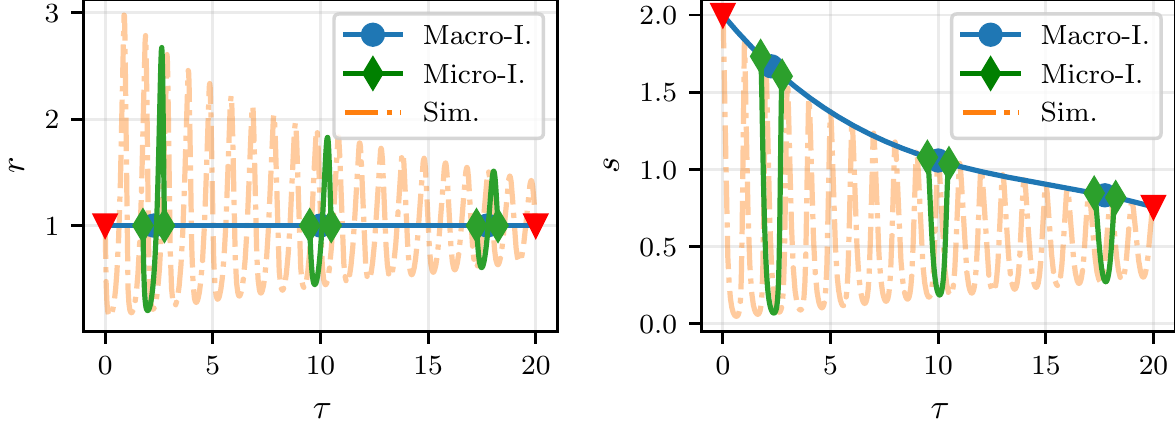}%
    \vspace{-7pt}%
    \caption{Found Solution to the OCP with Regularization \(\alpha = 2\cdot 10^{-5}\)}%
    \label{fig:PredatorPrey_OCP_Regularization} %
\end{figure}
\begin{figure}
    \centering%
    \vspace{4pt}
    \includegraphics*[width = \linewidth]{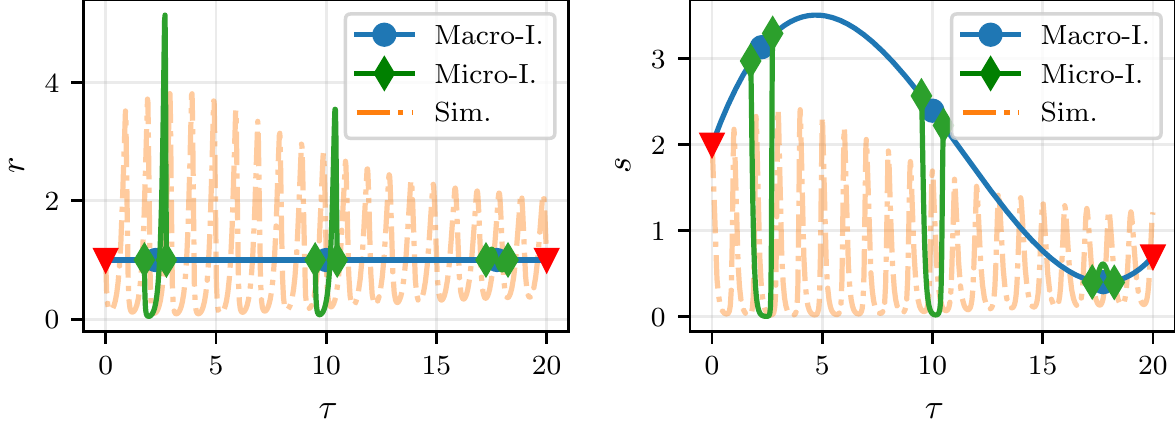}%
    \vspace{-7pt}%
    \caption{Found Solution to the OCP without Regularization \(\alpha = 0\)}%
    \label{fig:PredatorPrey_OCP_NoRegularization}%
\end{figure}
Fig. \ref{fig:PredatorPrey_OCP_Regularization} shows the optimal trajectory of the envelope in the states \(r\) and \(s\), the micro-integrations, as well as an a-posteriori performed simulation of the systems with the found optimal controls. 
With \(\alpha = 2\cdot 10^{-5}\) we observe that the trajectory and endpoint of the envelope matches the simulation very well.

In contrast, Fig. \ref{fig:PredatorPrey_OCP_NoRegularization} shows the found trajectories for \(\alpha = 0\).
The computed collocation trajectory reaches the target point \(\bar{x}_\mathrm{target}\), the tracking objective is thus zero, but the endpoint of the exact simulation \(x_{\mathrm{end},\mathrm{sim}}\) does not at all resemble the integration endpoint \(x_\mathrm{end}\).
This is because the optimizer chooses the optimization variables such that the envelope varies strongly.  
This then results in a large simulation error of the method, which then drives the integration endpoint to the target. 
The optimizer `used' the integration error to minimize the objective by driving the trajectories out of the allowed regions.
 
The choice of the regularization factor is important since it is a trade-off between smoother solutions that are not prone to the integration error and the original optimization objective. To investigate this, we solve the OCP for a range of \(\alpha \in [10^{-7},10^{-2}]\). 
As a baseline, we use the solution of the full OCP \eqref{eq:FullOCP} with the same polynomial control parametrization but without the regularization term whose optimal objective can be interpreted as the lower bound. 

\begin{figure}
    \centering
    \includegraphics*[width = \linewidth]{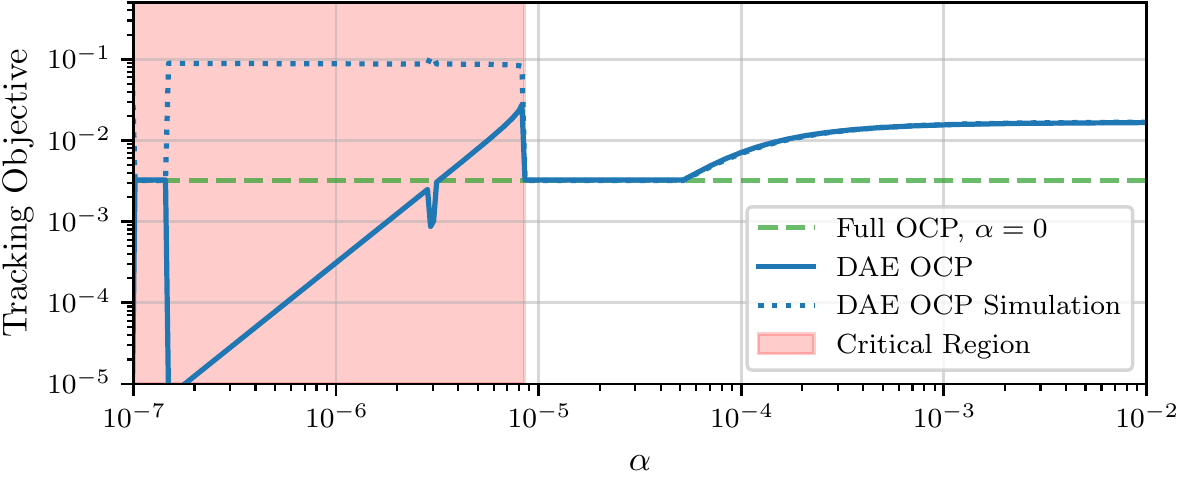}
    \vspace{-14pt}
    \caption{Tracking objective of the DAE OCP compared to an exact simulation result for the optimal controls for different \(\alpha\).}
    \label{fig:PredatorPrey_RegVsObj3}
\end{figure}
Fig. \ref{fig:PredatorPrey_RegVsObj3} shows the relationship between the choice of the regularization factor \(\alpha\) and the discrepancy between the `optimal' tracking objective and the simulated a-posteriori tracking objective. 

We observe that for large \(\alpha\) the regularization impacts the optimal trajectory considerably, such that the tracking objective is larger than necessary.
For decreasing \(\alpha\), the objective of the DAE OCP approaches the baseline objective; there is a `sweet-spot' region where the regularized and the full exact OCP coincide.
If \(\alpha\) is too small, the DAE OCP, through the erroneous integration explained above, converges to even lower objectives. In this region, the simulation does not match the found envelope trajectories. 

\section{Fuel-Optimal Orbit Transfers of a Low-Thrust Satellite}
\label{sec:Satellite}
As a second example, we consider a problem similar to the problem formulated by \cite{BETTS2000}. 
We want to find the optimal thrust direction and magnitude to transfer a satellite from a circular park orbit at \SI{500}{km} with an inclination of \(87.4\si{\degree}\) to a target mission orbit at \SI{1400}{km} with the same inclination but an eccentricity of \(0.0011\). 
Since the satellite has a very-low thrust engine, it requires a large number of orbits to reach the target orbit. 
The satellite has an initial mass of \(m_0 = \SI{160}{kg}\), a maximum thrust of \(T_0 = \SI{0.2}{\newton}\) and a specific impulse of \(I_\mathrm{sp} = \SI{1600}{s}\), as in \cite{BETTS2000}.
We model the satellite's state using modified equinoctial elements \(p,f,g,h,k,L\), where only the semi-late rectum \(p\) and the true longitude \(L\) have intuitive physical meaning\cite{Yang2009}. Their relationship to classical orbit elements is given in \cite{BETTS2000}.
Also including a state for the satellite mass \(m\) as well as a clockstate \(t\), we obtain the state vector and control vector
\begin{align}
    x = [p,f,g,h,k,L,m,t] \in \R^8 , && u \in \R^3.
\end{align}
The state dynamics are given by 
\begin{align}
    f(x,u) := \begin{bmatrix}
        A(x) \Delta(x,u) + b(x) \\
        - \frac{T_0}{g_\mathrm{e}I_\mathrm{sp}} \sqrt{u^\top u + \epsilon^2}\\
        1
    \end{bmatrix} \label{eq:SatelliteDynamics}%
\end{align}
where  \(g_e = \SI{9.81}{\metre\per\square\second}\) is the gravitational acceleration at sea level\cite{Yang2009}. Information on the matrices \(A\) and \(b\) can be found in \cite{BETTS2000}. The perturbation acceleration
\begin{align}
    \Delta(x,u) = \Delta_\mathrm{G}(x) + \frac{T_0}{m} u
\end{align}
contains a term for the gravitational perturbations \(\Delta_\mathrm{G}\) as described by \cite{BETTS2000}.
These perturbations are due to earth's oblateness, they dominate the satellite's thrust\cite{BETTS2000} and give rise to the highly oscillatory state trajectory.
The acceleration by the thruster is characterized by the thrust vector \(u = [u_1, u_2,u_3]^\top \in \R^3\) that satisfies
\begin{subequations}
    \begin{alignat}{2}
        -1 \leq & u_i &&\leq 1, \qquad i = 1,2,3, \label{eq:SatelliteControlBounds}\\
        & u^\top u &&\leq 1. \label{eq:SatelliteControlInequ}
    \end{alignat}
\end{subequations} %
For numerical reasons, we include a small parameter \(\epsilon=10^{-5}\) to smoothen the nonlinear square root in \eqref{eq:SatelliteDynamics}.
We constrain the classical elements at the endpoint of the trajectory to be within certain bounds of the target values, as explained in \cite{BETTS2000}.

We discretize the trajectory into orbits with the phase conditions
\begin{subequations}%
    \begin{alignat}{2}
        &L^\smallMinus &&= L_0\\
        &L^\smallPlus  &&= L_0 + 2\pi
    \end{alignat} \label{eq:SatelliteConditions} %
\end{subequations} %
on the phase state \(L\), where \(L_0\) is the initial value of \(L\) and \(q,b^\smallMinus,b^\smallPlus\) are chosen accordingly.
This type of conditions (B) is necessary since the dynamics \eqref{eq:SatelliteDynamics} are \(2\pi\)-periodic in the state \(L\).  

For the micro-integration of an orbit, we divide it into \(N_\mathrm{ctrl} = 30\) piecewise constant control intervals that each integrate the dynamics using one step of a 3-stage Gauss-Legendre (GL) collocation scheme of order 6, which results in \(U(\tau) \in \R^{3 \times 30}\). The algebraic variables now also contain the variables of the collocation scheme. 

We search for a fuel optimal transfer over \(578\) orbits, that we split into two intervals of \(N=289\) orbits each. 
In each interval, we integrate the DAE dynamics using a 5-stage GL collocation, and use a constant parametrization of the controls \(U(\tau)\) with \(d_u = 1\). 
The macro-integration variables are initialized by a linear interpolation between the initial state and target final state, the micro-integration variables by simulating a single orbit at the interpolation point. The control variables for all intervals are initialized as \([0,\ 0.1,\ 0]^\top\).

\begin{figure}
    \vspace{4pt}
    \centering
    \includegraphics*[width = \linewidth]{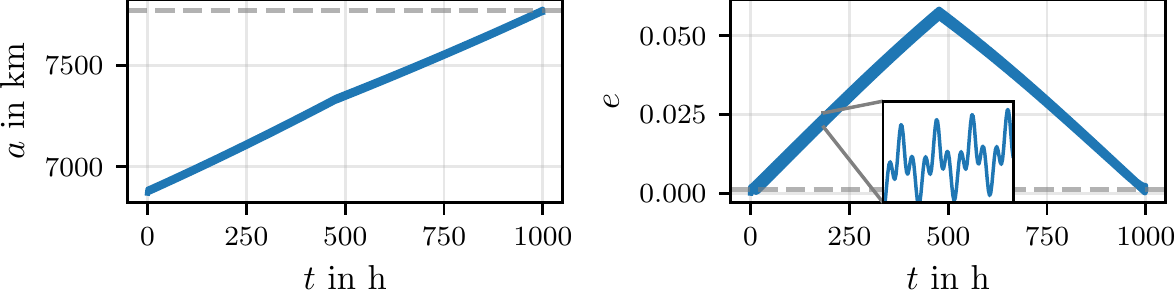}
    \vspace{-17pt}%
    \caption{Semi-major axis and eccentricity of a simulation of the found optimal controls with the target value indicated by the dashed line}%
    \label{fig:ClassicalElements}%
\end{figure}

\begin{figure}
    \vspace{4pt}
    \centering
    \includegraphics*[width = \linewidth]{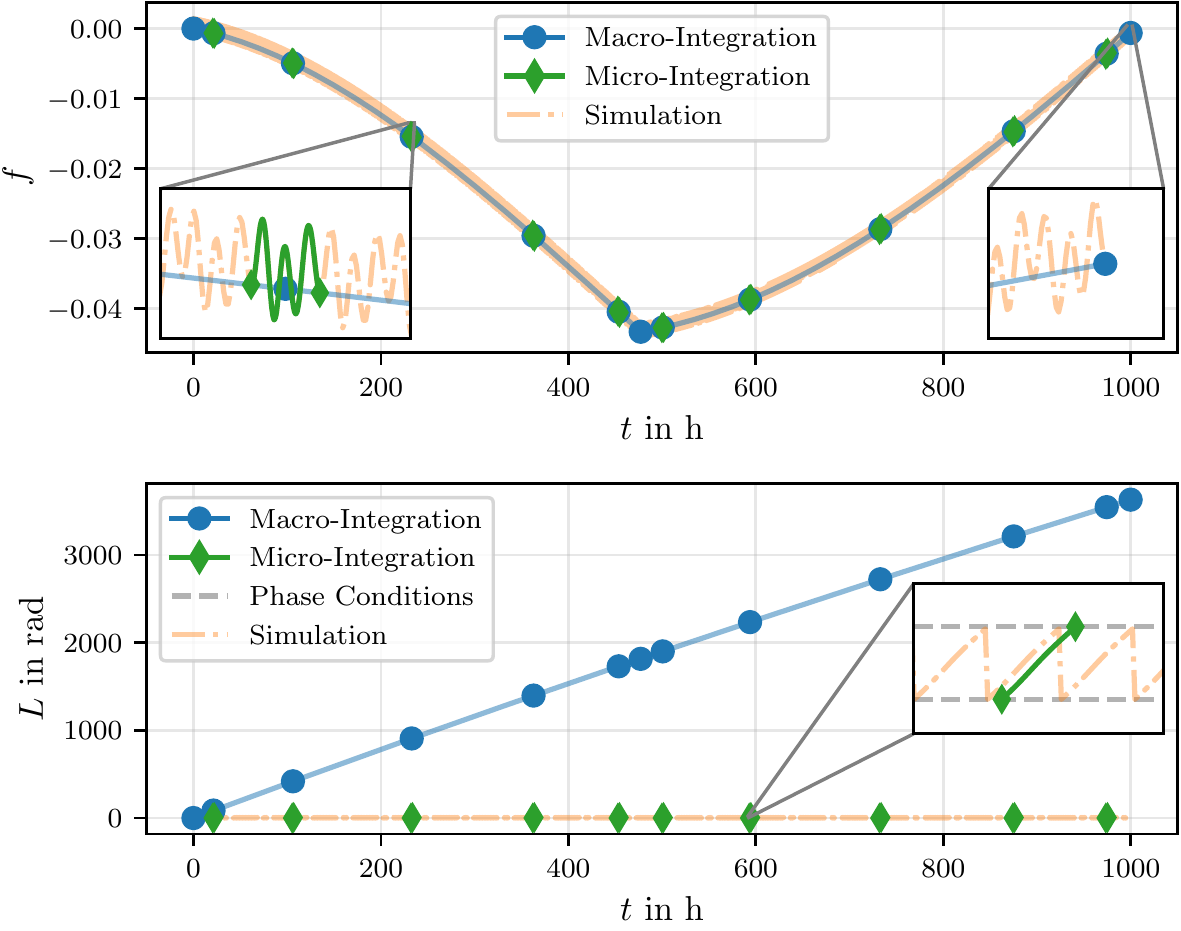}
    \vspace{-17pt}%
    \caption{The highly oscillatory trajectory of the state \(f\) and \(L\) as computed by the DAE method as well as a simulation of the found optimal controls. The state \(L\) is subject to the phase conditions \eqref{eq:SatelliteConditions}, in this state the macro-integration values are not tied to values of the start and endpoint of the micro-integration.}
    \label{fig:SatelliteSingleState}
\end{figure}

We solve the NLP with a total of \(17\ 334\) variables and  \(17\ 222\) constraints in about 20 seconds.
The results are presented in Fig. \ref{fig:ClassicalElements} and \ref{fig:SatelliteSingleState}. 
For reference, we simulate the found optimal controls with high accuracy.
We observe that the envelope agrees with the simulation over the orbits and at the endpoint, with a maximum relative error of \(8.85\cdot10^{-5}\) in the state \(g\).
The optimal final mass is \SI{155.48}{kg}, which almost coincides with the value found by \cite{BETTS2000}. 
Due to memory limitations, we were not able to solve the 32-times bigger full optimization problem that for the chosen control and state discretization has about \(555\ 000\) variables.

\section{Conclusion}
\label{sec:Conclusion}
We presented a novel method to treat highly oscillatory optimal control problems.
The state and control trajectories of the problem are discretized into whole cycles using phase conditions. 
Instead of simulating all cycles, we instead integrate a differential-algebraic equation that approximates the envelope dynamics of the state trajectory.
As a numerical example, we used the presented method to efficiently optimize the orbit transfer of a low-thrust satellite. 
Future work will adress a more rigourous analyis and theoretical justification of the method as well as the exploration of more use-cases. 
\bibliography{sources}

\end{document}